\documentclass[reqno]{amsart}
\usepackage{amssymb}
\usepackage{tikz-cd}
\usepackage[colorlinks, allcolors=blue]{hyperref}

\newtheorem{theorem}{Theorem}[section]

\newtheorem{lemma}[theorem]{Lemma}
\newtheorem{corollary}[theorem]{Corollary}

\newtheorem{definition}[theorem]{Definition}

\theoremstyle{remark}
\newtheorem{remark}[theorem]{Remark}

\numberwithin{equation}{section}

\begin{document}

\title[The basic representation of the DAHA at critical level]{On the basic representation of the double affine Hecke algebra at critical level}

\author{J.F.  van Diejen}

\address{
Instituto de Matem\'aticas, Universidad de Talca,
Casilla 747, Talca, Chile}

\email{diejen@inst-mat.utalca.cl}

\author{E. Emsiz}

\address{Delft Institute of Applied Mathematics,  Delft University of Technology, P.O. Box 5031, 2600 GA Delft, The Netherlands}
\email{e.emsiz@tudelft.nl}

\author{I.N. Zurri\'an}

\address{Departamento de Matem\'atica Aplicada II, Universidad de Sevilla.
EPS c/ Virgen de Africa 7, 41011 Sevilla, Spain}

\email{ignacio.zurrian@fulbrightmail.org}


\thanks{This work was supported in part by the {\em Fondo Nacional de Desarrollo
Cient\'{\i}fico y Tecnol\'ogico (FONDECYT)} Grant \# 1210015}

\date{September, 2022}

\begin{abstract}
We construct the basic representation of the double affine Hecke algebra at critical level $q=1$ associated to an irreducible reduced affine root system $R$ with a reduced gradient root system. 
For $R$ of untwisted type such a representation was studied by Oblomkov~\cite{obl:double} and further detailed by  Gehles~\cite{geh:properties} in the presence of minuscule weights.

\end{abstract}

\maketitle

\section{Introduction}


Our aim is to construct  a monomorphism of the double affine Hecke algebra at critical level  $\mathcal H$ associated  to any irreducible reduced affine root system $R$ with a reduced gradient   root system.
To this end we adapt corresponding results of Cherednik~\cite{che:double} and Macdonald 
\cite{mac:affine} (but see also \cite{hai:cherednik}) for general Macdonald parameter values $q$ (not equal to a root of unity) to the critical level $q=1$, with the aid of techniques from \cite[\text{Sec. 3}]{obl:double},
\cite[\text{Ch. 2.1}]{geh:properties} and \cite[\text{Ch. 4.3}]{mac:affine}.

For $R$ of untwisted type the basic representation at critical level was considered by 
Oblomkov~\cite[\text{Secs. 3 and 5}]{obl:double} in his study of the center and the ring-theoretical properties of $\mathcal{H}$; details of this construction were worked out in \cite[\text{Sec. 2.1}]{geh:properties} for the case of a nontrivial $\Omega$ (where $\Omega$ denotes the subgroup of the extended affine Weyl group consisting of elements of length zero).

The underlying idea in the later case is that the generator of $\mathcal H$ associated to the affine simple reflection can be suppressed in the defining relations because it can be expressed as a conjugation of a generator associated to a finite simple reflections by an element  in $\Omega$. 
In this note we use a different technique, which is based on the reduction to the subalgebras associated with the parabolic subgroups of the affine Weyl group
that fix a vertex of the Weyl alcove (see  Remark~\ref{rem:Omega-trivial} for further details).

One reason that   double affine Hecke algebras at critical level associated to affine root systems  and their representations are useful is that they can be used to compute affine Pieri rules for Hall-Littlewood polynomials and also govern the underlying symmetry structures of certain discretizations of  quantum integrable systems with delta-potentials on 
the root hyperplanes of an affine root system~\cite{ die-ems:discrete,die-ems-zur:completness}  (see  Remark~\ref{polynomial-rep} for more details).

The material is organized as follows. After recalling preliminaries concerning the definition of the double affine Hecke algebra at critical level in Section~\ref{daha}, the basic representation is constructed and its proof is provided in Section~\ref{sect:basic-rep}. 








\section{The double affine finite Hecke algebra at critical level}\label{daha}

\subsection{The affine root system}
Let $V$ be a real finite-dimensional Euclidean vector space with inner product $\langle \cdot ,\cdot \rangle$. We identify the set $\hat V$ of affine-linear functions on $V$ with $V+ \mathbb R$ via $(v+r)(x):=\langle v,x\rangle+r$. We also need the gradient map from $\hat V$ to $V$, defined by $(v+r)':=v$. 
We fix an irreducible reduced affine root system  $R\subset \hat V$ 
with a reduced gradient  root system $R_0:=R'=\{a'\mid a\in R\}\subset V$ 
of full rank in $V$. 
We write $Q$, $P$, and $W_0$, for the root lattice, the weight lattice, and the Weyl group associated with $R_0$. We also fix a  choice of positive roots $R_0^+$. The semigroup of the root lattice generated by $R_0^+$
is denoted by $Q^+$ whereas  $P^+$  stands for the corresponding cone of dominant weights.

We  fix a simple basis $\alpha_1,\ldots ,\alpha_n$ for $R_0^+$ and extend it to a  basis   $a_0,\ldots ,a_n$  of $R$, with  $a_0:=\alpha_0+ c$, $a_j:=\alpha_j$ ($j=1,\ldots ,n$) and $c>0$.   
Here  $\alpha_0=-\vartheta$  if  $R$ is \emph{twisted} and   $\alpha_0=-\varphi$ if $R$ is \emph{untwisted}, where $\varphi$ and $\vartheta$  denote the highest root and the highest short root of  $R_0^+$, respectively. If $R_0$ is simply-laced, then we will  call all roots long and short.
(So if $R_0$ is simply-laced, then $R$ is twisted and untwisted).
 The affine root system $R$ consists of the set  of elements $\alpha +  m_\alpha r c$ $(\alpha\in R_0, r\in \mathbb Z)$, where $m_\alpha=1$  if  $\alpha$ is short root and $m_\alpha=\frac{\langle\varphi,\varphi\rangle}{\langle \alpha_0,\alpha_0\rangle}$ if $\alpha$ is long. 

\subsection{The affine Weyl group}
An affine root $a=\alpha+m_\alpha rc\in R$ gives rise to an affine reflection $s_a:V\to V$ across the hyperplane $V_a:=\{ x\in V\mid a(x)=0\}$ given by 
$s_a(x):=x-a(x)\alpha^\vee$. The {\em affine Weyl group} $W^a$ is defined as the group generated by the affine reflections $s_a$, $a\in R$ and contains the finite Weyl group $W_0$ as the subgroup fixing the origin.
For any $y\in V$ we  denote by $t_y:V\to V$ the translation determined by the action $t_y(x):=x+y$.  
Let $\hat R_0 := \{\hat \alpha \mid \alpha\in R_0\}$ with $\hat \alpha:=m_\alpha \alpha^\vee$ (so $\hat R_0 =  \frac{2}{\langle \vartheta,\vartheta\rangle}R_0$ in the twisted case and $\hat R_0 =  R_0^\vee$ in the untwisted case, i.e. $\hat R_0=R_0^\vee$ for simply-laced $R_0$).
Then the elements of $W^a$ can be written as $vt_{ c\lambda}$ with $v\in W_0$ and $\lambda$ in $ \hat Q$, i.e. $W= W_0 \ltimes t(c{\hat Q})$, where $\hat Q$ denotes the root lattice of $\hat R_0$. 
By extending the lattice from 
${\hat Q}$ to the weight lattice ${\hat P}$  of $\hat R_0$ one arrives at the {\em extended affine Weyl group}
$ W=W_0\ltimes t(c{\hat P})$ consisting of group elements of the form $vt_{c\lambda}$ with $v\in W_0$ and $\lambda\in {\hat P}$. 
The action of $W$ on $V$ induces a dual action on the space $\mathcal{C}(V)$ of functions $f:V\to\mathbb{C}$ given by
$(wf)(x):=f(w^{-1}x)$ ($w\in W,\, x\in V$). 
 The affine root system $R\subset \mathcal{C}(V)$ is stable with respect to this dual action.

We denote by $\Omega$ the subgroup of elements of $ W$ that permutes the simple affine roots. So in particular $ua_j=a_{u_j}$ ($u\in \Omega$), and thus $us_ju^{-1}=s_{u_j}$,  where $j\to u_j$ encodes the corresponding permutation of the indices
 $j=0,\ldots ,n$. The upshot is that $W=\Omega\ltimes W^a$, with $W$ being normal in $W$,
 whence $\Omega$ is a finite abelian subgroup: $\Omega\cong  W/W^a \cong {\hat P}/{\hat Q}$.
 The extended affine Weyl group $W$ can now be presented as the group generated by the commuting
elements from
$\Omega$ and the simple reflections $s_0,\ldots s_n$, subject to the relations
\begin{subequations}
\begin{equation}\label{W-rel1}
u s_j u^{-1}=s_{u_j}    \quad u\in \Omega ,\ j=0,\ldots ,n,
\end{equation}
and
\begin{equation}\label{W-rel2}
(s_js_k)^{m_{jk}}=1 \quad  j,k =0,\ldots ,n.
\end{equation}
\end{subequations}
Here $m_{jk}=1$ if $j=k$ and
$m_{jk} \in \{   2,3,4,6 \}$  if $j\neq k$
(with the provision that for $n=1$ the order
$m_{10}=m_{01}=\infty$).
It follows that
any $w\in W$ can be decomposed (nonuniquely) as
\begin{equation}\label{red-exp}
w=u s_{j_1}\cdots s_{j_\ell},
\end{equation}
with $j_1,\ldots,j_\ell\in \{ 0,\ldots ,n\}$, and $u\in\Omega$.
The length $\ell (w)$ is defined as the minimum number of reflections $s_j$ ($j=0,1,\dots,n$) involved in any decomposition \eqref{red-exp} of $w$. Any decomposition \eqref{red-exp} with $\ell=\ell(w)$ is called a {reduced decomposition} of $w$. The group $\Omega$ can thus be also characterized as the group of all elements $ W$ that have length zero.

\subsection{The double affine Hecke algebra at critical level}

Let $\tau_\alpha$ ($\alpha\in R_0$) be a family of  indeterminates such that 
$\tau_{w\alpha}=\tau_{\alpha}$ for all $\alpha\in R_0$. 
We  denote by $\mathbb L$ the complex algebra of Laurent polynomials in the indeterminates $\tau_\alpha$ ($\alpha\in R_0$). It will be also useful to extend the definition of indeterminates by setting $\tau_{a}:=\tau_{a'}$ for $a\in R$ and  $\tau_j:=\tau_{\alpha_j}$ for $j=0,1,\dots,n$.

The {\em double affine Hecke algebra at critical level}  $\mathcal{H}$ associated to the affine root system $R$ is the unital associative $\mathbb{L}$-algebra  with invertible generators $T_u$ ($u\in\Omega$), $T_j$ ($j=0,\dots,n$) and the commuting elements $X^\lambda$ ($\lambda \in P$),
such that the following relations are satisfied

 \begin{subequations}\label{sub}
 \begin{equation}\label{Tu-rel}
 T_uT_{\tilde{u}}=T_{u\tilde{u}} \quad\text{and}\quad T_uT_j=T_{u_j}T_u  \qquad
 (u,\tilde{u}\in\Omega ,\, 0\leq j\leq n),
 \end{equation}
 \begin{equation}\label{quadratic-rel}
 (T_j-\tau_j)(T_j+\tau_j^{-1})=0\qquad (0\leq j\leq n),
 \end{equation}
 \begin{equation}\label{braid-rel}
 \underbrace{T_jT_kT_j\cdots}_{m_{jk}\ {\rm factors}}
 =\underbrace{T_kT_jT_k\cdots}_{m_{jk}\ {\rm factors}}\qquad (0\leq j\neq k\leq n),
 \end{equation}
 \begin{equation}\label{TuX-rel}
 T_u X^\lambda = X^{u^\prime\lambda} T_u \qquad (u\in \Omega,\, \lambda\in P),
 \end{equation}
 \begin{equation}  \label{TjX-rel}
 T_jX^\lambda = X^{s_j^\prime \lambda }T_j+(\tau_j-\tau_j^{-1})\frac{X^\lambda-X^{s'_j\lambda}}{1-X^{-\alpha_j}} \qquad
 (0\leq j\leq n,\,  \lambda\in P),
 \end{equation}
 \begin{equation}\label{X-rel}
 X^\lambda X^{\tilde{\lambda}} = X^{\lambda+\tilde{ \lambda}} \qquad ( \lambda,\tilde{\lambda}\in P ),
 \end{equation}
 \end{subequations}
where the number of factors $m_{jk}$ on both sides of the braid relation \eqref{braid-rel} is the same as the order of the corresponding braid relation \eqref{W-rel2} for $ W$. 
Here, for any element in $ W$ of the form $v t_{c\lambda}$ ($v\in W_0$, $\lambda\in \hat P$) we denote $(v t_{c\lambda})':=v$. Hence
 $s'_j=s_j$ for $j=1,\dots,n$ and $s'_0$ is the orthogonal reflection associated with $\alpha_0$.

\begin{remark}
Although the results in the note are formulated in terms of the algebra $\mathbb{L}$, all statements (and proofs) are also valid when the  indeterminants $\tau_\alpha$ are specialized to non-zero complex values. 
\end{remark}

We now turn to the construction of  the basic representation.


\section {The basic representation}\label{sect:basic-rep}


We denote the group algebra of $P$ by $ \mathbb{L}[P]$ and the corresponding natural basis elements (also) by $X^\lambda$, 
$\lambda \in P$. The finite Weyl group $W_0$ acts on $\mathbb{L}[P]$ via $wX^\lambda=X^{w\lambda}$. 
We will also need the localization $\mathcal A$ of $\mathbb{L}[P]$ by  the Weyl denominator  
$\delta:=\delta(R_0)=\prod_{\alpha\in R_0}(1-X^{\alpha}).$ Then the action of $W_0$ on $\mathbb{L}[P]$ extends to an action on  $\mathcal A$ via $\mathbb{L}$-isomorphisms. We write $f^w$ for the result of the action of $w\in W_0$ on $f\in\mathcal A$. 

\begin{definition}\label{smash-product}
The \emph{smash product}   $\mathcal A\, \# \, W$ is the associative unital $\mathbb L$-algebra  characterized by the following properties:
\begin{itemize}
\item[i)] $\mathcal A\, \# \, W$ contains $\mathcal A$  and the group algebra $\mathbb{L}[ W]$ as subalgebras,
\item[ii)] the multiplication map defines an isomorphism
$\mathcal A \otimes_{\mathbb L} \mathbb{L}[W] \to {\mathcal A}\, \# \, W$  of $\mathbb{L}$-modules, and
\item[iii)] we have the cross relations $(fv)(gw) = f g^{v'} vw$ for all $f,g\in {\mathcal A}$ and for all $v,w\in W$. 
\end{itemize}
\end{definition}

We will also need the following Demazure-Lusztig-type elements in  $\mathcal A\, \# \, W$:
\begin{equation}\label{dl}
\check T_j:=\tau_j s_j+ \frac{\tau_j-\tau_j^{-1}}{1-X^{-\alpha_j}} (1-s_j)  \qquad (j=0,\dots,n).
\end{equation}

\begin{theorem}[The Basic Representation of $\mathcal{H}$]\label{basic-rep}
The assignments $T_j \mapsto  \check{T}_j$ ($j\in\{0,\dots, n\}$), $T_u\mapsto u$ ($u\in\Omega$) and $X^\lambda \mapsto X^\lambda$ ($\lambda\in P$)  uniquely extend to an injective homomorphism $\jmath: \mathcal{H} \to \mathcal A\, \# \, W$ of $\mathbb{L}$-algebras. 
\end{theorem}

\begin{remark}\label{ore-extension}

The above monomorphism can  be lifted to an isomorphism after a suitable localization. For this  consider the extension  $\mathcal{H}_{\delta_\tau}$ of $\mathcal{H}$ obtained by adjoining the inverses of  $1-X^\alpha$ and $\tau_\alpha - \tau^{-1}_\alpha X^\alpha$ ($\alpha\in R_0$). In other words,   $\mathcal{H}_{\delta_\tau}$ is the Ore extension  obtained by  localizing  $\mathcal{H}$ at the $\tau$-deformed Weyl denominator 
$\delta_\tau= \prod_{\alpha\in R_0}(1-X^\alpha)(\tau_\alpha^{-1}-\tau_\alpha X^{\alpha}).$
Then the assignments in Theorem~\ref{basic-rep} uniquely extend
to an $\mathbb{L}$-isomorphism
$\jmath: \mathcal{H}_{\delta_\tau} \xrightarrow{\sim} \mathbb{L}[P]_{\delta_\tau} \, \# \, W.$
\end{remark}

\begin{remark}[The Polynonomial Representation of $\mathcal H$] \label{polynomial-rep} 
 For $c\in\frac{\langle\alpha_0,\alpha_0\rangle}{2}\mathbb{N}$ 
 the affine Weyl group $W$ acts on $\mathbb{L}[P]$ via $wX^\lambda=X^{w(\lambda)}$; this gives rise to a corresponding action of the Demazure-Lusztig operators~\eqref{dl} on $\mathbb{L}[P]$. By letting the monomials $X^\lambda\in\mathcal{H}$ ($\lambda\in P$) act  as multiplication operators the basic representation then reduces at such values for $c$ to a (unfaithful) representation of $\mathcal{H}$ on $\mathbb{L}[P]$.
In~\cite[Sect. 4]{die-ems:discrete} this polynomial representation  was, for $R$ of twisted type and the indeterminants $\tau_\alpha$ specialized to non-zero complex values, dualized so as to arrive at a representation of $\mathcal{H}$ given by translation operators and integral-reflection operators used to build quantum integrable systems in a Hilbert space of complex functions on the weight lattice $P$ (cf. also \cite[Sect. 4 and 5]{die-ems-zur:completness} for a corresponding construction associated with the double affine Hecke algebra of type $C_nC_n^\vee$ at critical level). These quantum integrable systems can be interpreted as lattice regularizations of quantum models with delta potentials related to the Cherednik algebra at critical level \cite{ems-opd-sto:trigonometric}.
\end{remark}
The rest of this section is devoted to the proof of  Theorem~\ref{basic-rep}. We start with the construction of the homomorphism.  For this we have to show that the relations \eqref{sub} holds when $T_j$ is replaced by   $\check{T}_j$ ($j=0,\dots, n$), $T_u$ by $u$ ($u\in \Omega$) and $X^\lambda\in \mathcal H$ by $X^\lambda \in\mathbb L[P]$ $(\lambda\in P$). 
For any affine root $a$ with gradient $\alpha$ let $\check T(a)=\tau_a s_a+ \frac{\tau_a-\tau_a^{-1}}{1-X^{-\alpha}} (1-s_a)$, so in particular $\check T_j = \check T(a_j)$. 
 Relations  \eqref{TuX-rel} and \eqref{X-rel} are trivial while relation \eqref{Tu-rel} is a consequence of $w \check T(a) w^{-1}= \check T(wa)$ ($w\in W$) and $u a_j = a_{u_j}$ ($j=0,\dots, n$).
Relation \eqref{TjX-rel} follows from a direct calculation.

We shall now concentrate on the quadratic relations  \eqref{quadratic-rel} and the braid relations \eqref{braid-rel}.
 For this we consider the affine Hecke algebra $H$ associated to $R$, i.e.  the unital associative $\mathbb{L}$-algebra  with invertible generators  $T_0, T_1,\dots, T_n $ subject to the relations Eqs. \eqref{quadratic-rel}  and \eqref{braid-rel}.  The relations  \eqref{quadratic-rel} and \eqref{braid-rel} are an immediate  consequence of the following lemma. 
\begin{lemma}\label{affine-homomorphism}
The assignments $T_j \mapsto  \check{T}_j$ ($j\in\{0,\dots, n\}$) uniquely extend to a homomorphism $\jmath: H \to \mathcal A\, \# \, W$ of $\mathbb{L}$-algebras. 
\end{lemma}

The idea of the proof of the above lemma is to reduce it to the   subalgebras associated with the parabolic subgroups of the affine Weyl group that fix a vertex of the Weyl alcove. These subalgebras are of finite type and this leads us to recall the definition of the finite Hecke algebra $H(R_0)$ associated to $R_0$ (in the rest of this paragraph we will also allow $R_0$ not to be irreducible):  it is  the unital associative $\mathbb{L}$-algebra  with invertible generators  $T_1, T_2,\dots, T_n $ subject to the relations Eqs. \eqref{quadratic-rel}  and \eqref{braid-rel} that do not involve $T_0$.   
We will also need the smash product $\mathcal A_Q \, \# \, W_0$, where $\mathcal A_Q= \mathbb{L}[Q]_{\delta(R_0)}$, and which is defined analogously as $\mathcal A\,\# \, W$ but with the  cross relations iii) in Definition~\ref{smash-product} replaced by $(fv)(gw) = f g^{v} vw$ for  $f,g\in \mathcal A_Q$ and $v,w\in W_0$.  
The elements  $\check{T}_j$ ($j\in\{1,\dots, n\}$) can be  naturally  interpreted as elements of  $ \mathcal A_Q \, \# \, W_0$.  
The proof of Lemma~\ref{affine-homomorphism} makes use of the following lemma:
\begin{lemma}\label{basic-rep-H0}
The assignments
$T_j \mapsto  \check{T}_j$ ($j=1,\dots, n$)  uniquely extend to a homomorphism $H(R_0) \to 
\mathcal A_Q\, \# \,W_0$ of $\mathbb{L}$-algebras. 
Moreover, this holds for any (not necessarily irreducible) reduced crystallographic root system $R_0$  of full rank in $V$. 
\end{lemma}

This lemma is an immediate consequence of the Bernstein-Lusztig relations \cite{lus:affine} in the affine Hecke algebra, see e.g.  \cite[(4.3.3)]{mac:affine}.  Although in loc. cit. it was assumed that $R_0$ was irreducible,  it actually holds for all finite crystallographic root systems $R_0$ of full rank: for this observe that if $R_0$  is decomposed into  disjoint, irreducible and orthogonal subsystems $\Sigma_1  \cup \dots \cup \Sigma_\ell$, then $H(R_0)$ decomposes as the tensor product of the $H(\Sigma_k)$, $k=1,\dots, \ell$  and similarly for  $\mathcal A_Q\, \# \,W_0$.

Let us now return to the proof of Lemma~\ref{affine-homomorphism}. For any $k\in\{0,1,\dots, n\}$ consider 
the parabolic $\mathbb{L}_k$-subalgebra $H_k$ of $H$ generated by $T_j$ $(j\neq k$), where $\mathbb{L}_k$ denotes the complex subalgebra of $\mathbb{L}$ generated by $\tau_j^{\pm 1}$ ($j\neq k$).
It suffices to show that, for any (fixed) $k$, the assignment $T_j\mapsto \check T_j$, $j\neq k$ extends to an homomorphism $H_k\to \mathcal A\,\# \, W$ of  $\mathbb{L}_k$-algebras.
For this we introduce the parabolic subgroup $W_k$ of $W$ generated by  $s_j$ ($j\neq k$) and set
 $R_k = \{w'\alpha_j \mid w\in W_k, j\neq k\}$. Then $R_k\subset R_0$ is  a finite root system of rank $n$ in $V$ (although  in general not  irreducible) with a basis of simple roots given by $\alpha_j$,  $j\neq k$.   The map $s_j \mapsto s'_j$, $j\neq k$ defines a Weyl group isomorphism $W_k\to W_0(R_k)$ and this isomorphism induces a natural homomorphism $\imath_k:H_k\to H(R_k)$ of $\mathbb{L}_k$-algebras defined by $T_j\mapsto T_j'$ where $T_j'$, $j\neq k$ denote the natural generators of $H(R_k)$.
 Let  $Q_k\subset Q$ be the root lattice of $R_k$ in $V$. If we denote the $\mathbb{L}_k$-subalgebra of $\mathcal A \, \# \, W$ generated by $\mathcal A_{Q_k}:=\mathbb{L}_k[Q_k]_{\delta(R_k)}$ and $W_k$
by $(\mathcal A \, \# \, W)_k$, then  $f w \mapsto f w'$ ($f\in A_{Q_k}$, $w\in W_k$) extends to an $\mathbb{L}_k$-isomorphism $\ell_k: (\mathcal A \, \# \, W)_k \to A_{Q_k}\, \# \, W_0(R_k)$.
By applying  Lemma \ref{basic-rep-H0} (with $R_0$ replaced by $R_k$) we obtain an algebra isomorphism $\jmath_k: H(R_k) \to \mathcal A_{Q_k}\,\#\, W_0(R_k)$. 
The proof of Lemma~\ref{affine-homomorphism} 
 now follows by observing that the homomorphism $\ell_k^{-1}\circ\jmath_k\circ\imath_k$ of $\mathbb{L}_k$-algebras maps $T_j$ to the Demazure-Lusztig-type elements $\check T_j$ \eqref{dl} for $j\neq k$. So in particular  $\jmath_k \circ \imath_k = \ell_k \circ \jmath_{| H_k}$. 
The situation can be visualized in the following commutative diagram:
\begin{equation}
\begin{tikzcd}
    H_k \arrow{d}{\imath_k} \arrow[dashed]{r}{\jmath_{| H_k}}
    &  (\mathcal A\,\#\, W)_k \arrow{d}{\ell_k} \\
    H(R_k) \arrow{r}{\jmath_k}
&\mathcal A_{Q_k}\,\#\, W_0(R_k)  \end{tikzcd}
\end{equation}

This proves Lemma~\ref{affine-homomorphism} and therefore also  the existence of the homomorphism $\jmath$.



We now turn to the injectivity of $\jmath$. 
For a $w\in W$ with a reduced expression $
w=u \,s_{j_1}\cdots s_{j_\ell},$ let 
$T_w:=T_u\,T_{j_1}\dots\,T_{j_\ell}.$
From the defining relations of $\mathcal H$ it follows  that the elements $X^\mu  T_w$ (or alternatively $T_w X^\mu$) span $\mathcal H$ over $\mathbb L$. 
To show injectivey it therefore  suffices to show that the elements $X^\mu  \check{T}_w$ (or alternatively $ \check{T}_w X^\mu$), where $\check{T}_w:= \jmath(T_w)$
with $\mu\in P$ and $w\in  W$, are linearly independent in $\mathcal A\, \# \, W$ over $\mathbb{L}$.  
This can be proven  by mimicking Macdonald's proof  in \cite[(4.3.11)]{mac:affine}  (but see also \cite[Prop. A.4]{die-ems-zur:completness}) of a corresponding statement for $q\neq 1$.  It is essentially based on the fact that 
for any   $w\in  W$ one has
$\check{T}_w  = \sum_{v\le w} f_{vw}(X) v$, 
where  the coefficients $f_{vw}(X)$ belong to $ \mathcal{A}$ and that the leading coefficient $f_{ww}(X)$ is non-zero (and where $\le$ refers to the Bruhat partial order on $W$, i.e. $v\le w$ iff a reduced expression for $v$ can be obtained from a reduced expression for $w$ by deleting simple reflections \cite{bou:groupes, hum:reflection}).
Furthermore, this triangularity property  can be inverted to yield the decomposition
$w =\sum_{v\le w} \tilde f_{vw}(X) \check{T}_v  \in  \mathbb{L}[P]_{\delta_\tau} \, \# \, W$, 
where  the coefficients $\tilde f_{vw}(X)\in \mathbb{L}[P]_{\delta_r}$ and with non-zero leading coefficient  $\tilde f_{ww}(X) = f_{ww}(X)^{-1}$. This  implies that the homomorphism from Remark~\ref{ore-extension} is an isomorphism.

\begin{remark}\label{rem:Omega-trivial}
If the finite group $\Omega$  is non-trivial, then there is an alternative proof of the existence of the homomormphism $\jmath$.
It is based on the fact that the generator $T_0$ and the relations~\eqref{Tu-rel}-\eqref{TjX-rel} involving $T_0$ are \emph{redundant} in this case because  there exist an $u\in \Omega$ and  $j\in \{1,2,\dots, n\}$ such that $s_0=u s_j u^{-1}$, and therefore also  $ T_0=T_{u} T_j T_{u}^{-1}$ and  $ \check T_0=u \check T_j u^{-1}$. 
Since the quadratic relations  \eqref{quadratic-rel} and braid relations \eqref{braid-rel} not involving $T_0$ follow  from Lemma~\ref{basic-rep-H0}, the existence of the homomorphism $\jmath$ follows in a straightforward way. 
This argument was employed by Oblomkov~\cite{obl:double}  and Gehles~\cite{geh:properties} to arrive at Theorem~\ref{basic-rep} in the case that $R$ is of untwisted type.
 In~\cite{die-ems-zur:completness}  it was shown for the five-parameter double affine Hecke algebra at critical level of type $C_nC_n^\vee$ that the proof in question can be adapted to the case of a trivial $\Omega$ via a conjugation/translation trick borrowed from  \cite{ems-opd-sto:trigonometric} (where it was used to construct integral-reflection representations of the trigonometric Cherednik algebra at critical level). The idea  of the present proof should be viewed as a uniform technique that works for all affine root systems $R$, irrespective of whether $\Omega$ is trivial or not. 
\end{remark}

As a consequence of (the proof of) Theorem~\ref{basic-rep} we obtain as an immediate corollary the following PBW property of Cherednik  for $\mathcal H$:
\begin{corollary}
The elements $X^\lambda T_w$ (or alternatively $T_w X^\lambda$), with $\lambda\in P$ and $w\in  W$, form an $\mathbb{L}$-basis for $\mathcal{H}$. 
\end{corollary}
Together with the  Bernstein-Zelevinsky relations in the affine Hecke algebra  \cite[Prop. 3.7]{lus:affine} this corollary also yields the 
 PBW property  $\mathcal H\simeq \mathbb{L}[\hat P]\otimes_{\mathbb{L}} H(R_0)\otimes_{\mathbb{L}} \mathbb{L}[P]$ (isomorphism as $\mathbb{L}$-modules), cf.~\cite[Thm. 3.2.1 (p. 310)]{che:double}. 




\vspace{3ex}
\noindent {\bf Acknowledgments.} 
This work was supported in part by the \emph{Fondo Nacional de Desarrollo Cient\'iico
y Tecnol\'ogico} (FONDECYT) Grant \# 1210015 and VI PPIT-US. 
We thank the referee for the helpful constructive remarks and suggesting some improvements concerning the presentation.

\bibliographystyle{amsplain}

\end{document}